\documentclass[11pt]{amsart} 
\usepackage{amssymb}
\usepackage{natbib}
\usepackage[T1]{fontenc} 
\usepackage[utf8]{inputenc} 
\usepackage{mathtools}
\usepackage[unicode,colorlinks,
	citecolor=blue,
	linkcolor=blue,
	urlcolor=black,
]{hyperref} 

\usepackage{booktabs}

\usepackage{enumitem} 

\usepackage{xcolor,soul}

\definecolor{Dark Orchid}{RGB}{153,50,204}
\definecolor{Goblin Green}{RGB}{99,181,33}
\definecolor{Crayola Bittersweet}{RGB}{254,111,94}

\theoremstyle{plain}
\newtheorem{thm}{Theorem}[section]

\newtheorem{lem}[thm]{Lemma}
\newtheorem{cor}[thm]{Corollary}

\theoremstyle{definition}
\newtheorem{defn}[thm]{Definition}

\newtheorem{query}[thm]{Question}

\theoremstyle{remark}

\numberwithin{equation}{section}

\theoremstyle{remark}
\newcommand{\thistheoremname}{}
\newtheorem*{genremark*}{\thistheoremname}

\renewcommand{\epsilon}{\varepsilon}
\renewcommand{\phi}{\varphi}

\renewcommand{\succeq}{\succcurlyeq}

\DeclareMathOperator{\uh}{\upharpoonright}

\DeclareMathOperator{\ran}{ran}

\begin{document}

\title{The power of enumerations}

\author{Karthik Ravishankar}
\address{Department of Mathematics, University of Wisconsin--Madison, 480 Lincoln Dr., Madison, WI 53706, USA}
\email{kravishanka3@wisc.edu}
\urladdr{https://karthik-ravishankar-math.github.io}

\date{\today}

\makeatletter
\@namedef{subjclassname@2020}{\textup{2020} Mathematics Subject Classification}
\makeatother

\begin{abstract}
 We study the computational strength of enumerations of generic and random sets, particularly its interplay with computing these `typical' sets. We show that, while enumerations of $2$-generics can compute $2$-generics, this fails at higher levels---every $3$-generic has an enumeration which computes no weak $3$-generic. The analogous statement holds for randoms---every weak $2$-random has an enumeration which computes no $1$-random.
 
 We then analyze the total degrees above a quasiminimal enumeration degree, showing certain restrictions as well as certain positive results on the possible Turing degrees in such a cone. This answers a question of Slaman and Sorbi, by showing that there are quasiminimals with no minimal Turing degrees above them. We also show that while every total degree above some quasiminimals compute dimension $1$ sets, the total degrees above any quasiminimal fail to be a subset of the PA degrees.
\end{abstract}

\maketitle

\section{Introduction}
Genericity and randomness are notions of typicality. In some ways, being typical, they cannot encode much information. On the other hand, they do allow for certain weak forms of computational power and consequently require some power to compute. In the first half of this paper, we analyze whether certain notions of computational strength at or weaker than enumerability allow for computing generics or randoms. 

In Section~2, we analyze the relationship between three notions of computational strength---Escaping power, computing infinite subsets and enumerating a set alongside various levels of genericity. We show that these notions line up for some levels of genericity while they fail to do so for others, concluding with Theorem~\ref{thm:enumnosuf}, where we show that enumerating generics does not suffice to compute a generic at the level of weak $3$-generics.

In Section 3, we analyze the relationship between enumerating and computing a $1$-random. Lemma~\ref{lem:w2r} shows that enumerating a sufficiently random does not suffice to compute a $1$-random. However enumerating $1$-randoms does allow for computing effective Hausdorff dimension $1$ sequences.

Enumeration reducibility is a positive analog to Turing reducibility, where $X \leq_e Y$ if there is a c.e.\ set $\Gamma$ such that $x \in X$ if and only if $\exists \langle x,F\rangle \in \Gamma$ with $F \subseteq Y$---where $F$ is a finite set provided via its canonical index. This reducibility induces its corresponding degree structure known as the enumeration degrees $(\mathcal D_e, \leq_e )$. 

The Turing degrees embed into the enumeration degrees under the map $\iota : X \to X \oplus \overline X$. This map is degree theoretic and respects meets, join and jump. The degrees in the range of $\iota$ in $\mathcal D_e$ are called the total degrees, and they correspond to an isomorphic copy of $(\mathcal D_T, \leq_T )$ within $(\mathcal D_e,\leq_e)$. 

Every total degree above some set $X$ will contain an enumeration of $X$ (surjective map $f: \omega \to X$) in it. The total degrees above an enumeration degree characterize it as a consequence of Selman's theorem: $X \leq_e Y $ if and only if every enumeration of $Y$ computes an enumeration of $X$. Therefore the following class of degrees is of interest.
\begin{defn}\label{def:largeecone}
$X\subseteq \omega$ is called quasiminimal if $X$ is not c.e.\ and $Y\oplus \overline Y \leq_e X$ then  $Y$ is computable. Let $\mathcal C(X) = \{ A \subseteq \omega : X \text{ is c.e.\ in } A\}$ be the e-cone above $X$. We call such $\mathcal C(X)$ \emph{large} if $X$ is quasminimal. 
\end{defn}
Note that $C(X)$ is Turing invariant, since `being c.e.\ in' uses both positive and negative information of the oracle. Therefore $C(X)$ can be thought of as the collection of Turing degrees which enumerate $X$, which we call an $e$-cone in the Turing degrees. It is easy to see that $\mathcal C(X)$ is large if and only if $\mathcal C(X)$ is not contained in a non trivial Turing cone. 

Every $1$-generic is quasiminimal as is every weakly $2$-random \cite{randomquasi}. So the class of quasiminimal e-degrees are large in the sense of measure and category. Several other constructions such as halves of $\mathcal K$-pairs, function generics etc also result in quasiminimals. 

Slaman and Sorbi \cite{minimalabovequasi} construct a large e-cone with continuum many minimal Turing degrees above it. They then ask if every large e-cone contains a minimal Turing degree. Note that minimal degrees are certain initial segments of $\mathcal D_T$. Using the analysis of generics in Section~4, we answer a stronger question.
\begin{thm}\label{thm:initalsegment}
There is a large e-cone which does not contain any initial segment (excluding $0$) of the Turing degrees.
\end{thm}

The following theorem summarizes the main results in this paper about quasiminimal degrees and the total degrees above them.
\begin{thm}\label{thm:quasi}
Fix an $X\subseteq \omega$: 
\begin{enumerate}
\item There is a large e-cone contained in the upwards closure of $X$-generics.
\item There is a large e-cone contained in degrees of dimension $1$ relative to $X$ sequences.
\item There is a large e-cone which contains no 1-generic set.
\item There is a large e-cone which contains no weak $2$-random set.
\end{enumerate}
Moreover no large e-cone is contained in the PA degrees.
\end{thm}

\section{Levels of genericity}
In this section we study different levels of genericity and whether escaping power, computing infinite subsets or enumerating a generic can characterize when one can compute a generic. 

Note that if $X$ can enumerate $Y$, then it can compute an infinite subset of $Y$, and if it can compute an infinite subset of $Y$ it can compute an escaping (in fact dominating) function to $p_Y$---the principal function of $Y$. Given a set of strings $W$, we say $X \in 2^\omega$ meets $W$ if $\exists \sigma \prec X$ such that $\sigma \in W$. We say $X$ avoids $W$ if $\exists \sigma \prec X$ such that$\forall \tau \succ \sigma ( \tau \not \in W)$.
\begin{defn}  $X$ is weak $n$-generic if it meets every dense $\Sigma^0_n$ set of strings. $X$ is $n$-generic if it meets or avoids every $\Sigma^0_n$ set of strings. $X$ is $Y$-generic if it meets or avoids every $\Sigma_1^Y$ set of strings.
\end{defn}

It is easy to see that $n$-generic implies weak $n$-generic which in turn implies $(n-1)$-generic for every $n\geq 1$. Moreover observe that being $n$-generic is the same as being $0^{(n-1)}$-generic. For ease of notation, let us call a weak $n$-generic a $(n-\frac{1}{2})$-generic, so for example a weak $3$-generic would be a $2.5$-generic.
\begin{defn}
Let $m \in \{\frac{n}{2}:n\in \omega\}$. 

\begin{enumerate}
\item We say escaping power is \emph{sufficient} for $m$-genericity if there is a family $\mathcal F$ of functions such that if $f$ escapes all members of $\mathcal F$, then $f$ computes a $m$-generic. We say escaping \emph{characterizes} $m$-genericity if moreover every $m$-generic has its principal function escaping $\mathcal F$.
\item We say computing infinite subsets is sufficient for $m$-genericity if there is an $n$ such that any infinite subset of a $n$-generic computes a $m$-generic. We say computing infinite subsets characterizes $m$-genericity if infinite subsets of $m$-generics compute (possibly different) $m$-generics.
\item We say enumerations are sufficient for $m$-genericity if there is an $n$ such that every enumeration of an $n$-generic computes an $m$-generic. We say enumerations characterize $m$-genericity if every enumeration of an $m$-generic computes a (possibly different) $m$-generic.
\end{enumerate}
\end{defn}
The table below summarizes the relationships between levels of genericity and computational strength. 

The following lemma by Jockush shows that computing infinite subsets of a generic is equivalent to the escaping power of a generic for all integral levels of genericity. The unrelativized proof appears in \cite{escinf}, we include a relativized proof here for completeness. 

\begin{lem} \label{lem:infesc} (Jockush) $A$ computes an infinite subset of an $X$-generic if and only if $\exists f \leq_T A$ such that $f$ escapes all $X$ computable functions.
\end{lem}
\begin{proof}
For the forward direction, observe that the principal function of an $X$-generic will escape all $X$ computable functions as it is dense to escape them. Therefore so will the principal function of a subset. 

For the converse, let $f \leq_T A$ be strictly increasing and escape all $X$ computable functions. Let $B$ be such that $p_B = f$ i.e.\ $f$ is the principal function of $B$. We will build a $G \supseteq B$ such that $G$ is $X$-generic.

We build $G$ in stages as $G = \bigcup_s \sigma_s$. Let $\sigma_0=\epsilon$. 

At stage $s+1$, we have constructed $\sigma_s$ and want to extend this so that we meet or avoid $W_s^X$. Search for an extension $\tau\succ \sigma_s$ which is compatible with $B$ ($\forall n < |\tau| B(n)=1 \implies \tau(n)= 1$) and such that i)$\tau \in W_s^X$ or ii) for every $\gamma \succ \tau$ we have $\gamma \not \in W_s^X$. Set $\sigma_{s+1} = \tau$ for the first such $\tau$ found. This ends the construction.

Suppose we get stuck at stage $s+1$. Then for every $B$-compatible extension $\tau \succ\sigma_s$, and for every $\gamma \succ \tau$ with $\gamma \in W_s^X$, $\gamma$ must be $B$ incompatible, and moreover for every such $\tau$, atleast one such $\gamma$ exists in $W_s^X$. But then, there is an $n >|\tau|$ such that $\gamma(n)=0$ and $B(n)=1$. This will allow us build a function $h\leq_T X$ which dominates $f$: $h(0) = |\gamma|$ where $\gamma \succ \sigma_s, \gamma \in W_s^X$, and $h(n+1) = |\gamma|$ for the first $\gamma \in W_s^X$ such that $\gamma \succ \sigma_s^\frown 1^{h(n)}$. 
\end{proof}
Therefore columns $1$ and $2$ below are equivalent, when we restrict ourselves to the non-weak levels of genericity.
\begin{table}[ht]
\centering
\renewcommand{\arraystretch}{1.25}
\begin{tabular}{|c|c|c|c|}
\hline
 & Escaping & Infinite subset & Enumeration \\ \hline
Weakly 1-generic & Characterizes & Characterizes & Characterizes \\ \hline
1-generic        & Sufficient    & Sufficient    & Sufficient    \\ \hline
Weakly 2-generic & Characterizes & Characterizes & Characterizes \\ \hline
2-generic        & Sufficient    & Sufficient    & Sufficient    \\ \hline
Weakly 3-generic & No & No & No \\ \hline
\end{tabular}
\end{table}

We now analyze the different entries in the table.
\begin{enumerate}
\item The principal function of a weak $n$-generic escapes all $\Delta^0_n$ functions since it is dense to escape them. The converse is easy construction for $n=1$: If $f$ escapes all total computable functions, it computes a weak $1$-generic. This is extended to the case of $n=2$ in \cite{bhhi}, and so for weak $1$-generics and for weak $2$-generics, we get an equivalence for all three notions of computational strength. This then characterizes degrees which compute the corresponding generics and completely explains rows $1$ and $3$. 

\item By Lemma~\ref{lem:infesc}, columns 1 and 2 grant the same power, and so will always line up.
\item Not all non c.e.\ $\Delta^0_2$ sets compute a $1$-generic (for example, consider a minimal degree below $0'$). Therefore enumerating a $1$-generic does not characterize the degrees which compute a $1$-generic, and neither do the weaker notions of escaping or infinite subset. This explains row $2$.

\item In \cite{bhhi}, the authors show that escaping $\Delta^0_3$ functions is sufficient to compute a $2$-generic. By Lemma~\ref{lem:enumsnochar} this is not a characterization, fully explaining row $4$. 

\item  In \cite{bhhi}, it is shown that no amount of escaping power is sufficient to compute a weak $3$-generic. A modification of this argument can be used to show that every weak $3$-generic has an infinite subset which does not compute any weak $3$-generic. This is given below in Lemma~\ref{lem:subsetno3}. Finally Theorem~\ref{thm:enumnosuf} shows that every $3$-generic has an enumeration which computes no weak $3$-generic. This fully explains the last row.
\end{enumerate}

\begin{lem}\label{lem:enumsnochar}
For every $n\geq 1$, there is a set $X$ which enumerates an $n$-generic but does not compute an $n$-generic.
\end{lem}
\begin{proof}
By relativizing the construction of a $\Delta^0_2$ minimal degree, we see that every $0^{(n-1)}$ has a minimal cover $M_n$, which is $\Delta^0_2(0^{(n-1)})$.  The authors in \cite{delta2enumerates} show that every non-c.e.\ $\Delta^0_2$ set enumerates a $1$-generic set. By relativizing this argument, we get that every set which is $\Delta^0_2(0^{(n-1)}) - \Sigma^0_1(0^{(n-1)})$ and above $0^{(n-1)}$ enumerates an $n$-generic.  

Therefore $M_n$, being non-c.e.\ relative to $0^{(n-1)}$, can enumerate a $n$-generic. If $X \leq_T M_n$ was $n$-generic, then setting $X= X_0\oplus X_1$, we have by Van Lambalgen's theorem for genericity, $X_0$ is $n$-generic relative to $X_1$ and hence $1$-generic relative to $X_1^{(n-1)}\geq_T X_1\oplus 0^{(n-1)} \equiv_T M$, a contradiction to the fact that $X_0 \leq_T M$. Hence $M_n$ cannot compute a $n$-generic.
\end{proof}

Every set $X$ admits an encoding $A\subseteq \omega$ such that every infinite subset of $A$ computes $X$. For example, let $A$ consist of natural number codes for all finite initial segments of $X$. However, when we restrict ourselves to infinite subsets of generic sets, we can avoid encoding information. The argument below is an extension of the argument in \cite{bhhi} which shows that escaping power does not ensure computing a weak $3$-generic.
\begin{lem}\label{lem:subsetno3} Every weak $3$-generic has an infinite subset which computes no weak $3$-generic.
\end{lem}
\begin{proof}
We build a $0''$ computable function tree $f: \omega^{<\omega} \to \omega^{<\omega}$ such that none of the paths through $f$ compute a weak $3$-generic. A path through $f$ will be treated as an enumeration of a set, and let $str(f(\tau))\in 2^{<\omega}$ be the binary string of length $max(f(\tau))$ corresponding to the set $\{f(\tau)(n): n<|f(\tau)|\}$. The construction will ensure that every weak $3$-generic is a superset of a path through $f$. 

During the construction, we build dense $\Sigma_3^0$ sets of strings $\{U_e\}_{e\in \omega}$ and ensure that no prefix of $\varphi_e^{f(\tau)}$ is in $U_e$.  Call an extension $\gamma \succ f(\tau)$ \emph{locking} if $\gamma$ enumerates no element below the largest element enumerated by $f(\tau)$, and we say $\gamma$ has locking level $n$ if the new elements it enumerates are all $>n$. Once we define $f(\tau)$, we associate a c.e.\ set of possible extensions $V_{\tau,s}$ to it, starting with $V_{\tau,0} = \{f(\tau)\frown m: m> \text{max}(f(\tau))\}$.

Let $\{\tau_s\}_{s\in\omega}$ enumerate $\omega^{<\omega}$, such that if $\tau_s \prec \tau_t $ then $s<t$. We now describe the construction of $f$. At stage $s=0$, let $f = \emptyset$. At stage $s+1$, we assume $f(\tau_{t})$ is defined for all $t<s$. We define $f(\tau_s)$ as follows. Let $s+1 = \langle e,i\rangle$ and let $k =i+s+1$. 

We ensure that $U_e$ dense by adding an extension of length $k$ of every string in $2^i$ as follows: For each $t<s$, check if  $\exists \rho_t \in 2^k\exists^\infty \gamma \in V_{\tau_t,s}$ such that each $\gamma$ has a locking extension $\hat \gamma \succ \gamma$ for which $\varphi_e^{\hat \gamma} \uh_k = \rho_t$. If yes, we set $V_{\tau,s+1}$ to be the set of those extensions. Otherwise $\rho_t$ is undefined and we set $V_{\tau,s+1} = V_{\tau,s} - $ finitely many $\gamma$'s which have locking extensions $\hat \gamma\succ \gamma$ for which $\varphi_e^{\hat \gamma}\uh_k \downarrow$. We extend $U_e$ by all strings at level $2^k$ except the $\rho_t$ for $t<s$. 

Note that we forbid at most $s$ many strings, while each string at length $i$ has $2^{s+1}$ many extensions at length $k$, and so every string of length $i$ has an extension in $U_e$. Then set $f(\tau_s) = $ first element entering $V_{\tau_t,s+1}$, where $\tau_t\frown i = \tau_s$, and remove $f(\tau_s)$ from $ V_{\tau_t,s+1}$. 

Note that if $g \in [f]$, then $\varphi_e^g$ is not weak $3$-generic since it is either partial or it does not meet $U_e$ by construction. Therefore $g$ does not compute a weak $3$-generic. Now let $G$ be a weak $3$-generic, and suppose $str(f(\tau))\subseteq \sigma \prec G$. Consider the $\Sigma^0_3$ set $\{\gamma \in 2^{<\omega}:  \gamma \supseteq str(f(\tau\frown i)), i\in\omega \}$. Since we left plenty of gaps by locking our enumerations on infinitely many inputs, this is dense above $\sigma$. Therefore $G$ must meet it, and $G$ has an infinite subset enumerated by a path through $f$.
\end{proof}

We next build a tree of enumerations, where the significant added difficulty is that we have to make sure to enumerate \emph{all} the elements of a generic. The tree we build will have some paths which fail and may end up computing weak $3$-generics, but we will be able to argue that for $3$-generics, we can always find an enumeration through our tree which succeeds. 

\begin{thm}\label{thm:enumnosuf}
Every $3$-generic set $G$ has an enumeration that computes no
weak $3$-generic set.
\end{thm}

\begin{proof}
As in Lemma~\ref{lem:subsetno3}, we build a $0''$ computable function tree $f: \omega^{<\omega} \to \omega^{<\omega}$ of enumerations and for each $\tau \in \omega^{<\omega}$ we associate a c.e.\ set $V_{\tau,s} \subseteq \omega^{<\omega}$ of possible locking extensions. We will further ensure that every finite extension of $str(f(\tau))$ has a valid extension at the next level.

We simultaneously try to build dense $\Sigma_3^0$ sets of strings $U_{e}^\tau$, so that $U_{e}^\tau$ witnesses failure of weak $3$-genericity of $\varphi_e^g$ for every path $g \in [f]$ above $f(\tau)$. We will argue that, while all paths may not succeed, $3$-genericity will be sufficient to ensure success.

Let $\{\tau_s\}_{s\in\omega}$ enumerate $\omega^{<\omega}$, such that if $\tau_s \prec \tau_t $ then $s<t$. We now describe the construction of $f$. At stage $s$, let $f(\epsilon) = \epsilon$. At stage $s+1$ we proceed as follows.

Let $s+1 = \langle e,i\rangle$. Let $u$ be the least number such that $f(\tau_u)$ is undefined, and $I_s = \{i: f_s(\tau_i)\downarrow\}$. Let $k = i+ |I| +1$. For each $t\in I$, if $U_{e}^{\tau_t}$ has not been stopped, we work towards making it dense: check if $\exists \rho_t \in 2^k \exists^\infty \gamma \in V_{\tau_t,s}$ with locking extensions $\hat \gamma \succ \gamma$ for which $\varphi_e^{\hat \gamma}\uh_k = \rho_t$. If such a $\rho_t$ exists, we set $V_{\tau_t,s+1}$ to be an infinite set of such $\hat \gamma$'s (at most one per $\gamma \in V_{\tau_t,s}$). Note that in this case, we are free to `backfill' a locking extension with smaller numbers as needed since that will not affect convergence.

Otherwise $\rho_t$ is undefined and we define $V_{\tau_t,s+1}$ as follows: Mark each of the cofinitely many  $\gamma \in V_{\tau_t,s}$ which have no locking extensions which give us a convergence on $k$ bits via $\varphi_e$ as \emph{$e$-divergent} with locking level $max(\ran(\gamma))$, and add them to $V_{\tau_t,s+1}$. 
We also stop building $U_{e}^{\tau_t}$ for the rest of the construction. 

Extend each unstopped $U_{e}^{\tau_t}$ by all the strings of level $k$ except $\{\rho_i\}_{\{i:\tau_i\succeq \tau_t\}}$.
If $\tau_w$ is  the parent of $\tau_u$, we now extend $f$ as follows.

Let $\gamma$ be the least element in the pool $V_{\tau_w,s+1}$. Let $\{F_i\}_{i<m}$ enumerate all finite subsets $F\subseteq (max(\ran(f(\tau_w)),max(\ran(\gamma)))$, with $F_0 = \emptyset$. To each $F_i$, assign an enumeration string $\gamma_i \succ \gamma$ which precisely enumerates $F_i$ such that $\{\gamma_i\}_{i< m}$ are incomparable. We then choose $m$ strings $\{\tau_w^\frown j_i\}_{i\leq m}$ (including $\tau_u$) on which $f$ is undefined, and set $f(\tau_i) = \gamma_i$. 

Moreover, if $\gamma_i$ as an extension of $\gamma$ has locking level $n$, it inherits the $e$-divergement labeling of $\gamma$ with locking level $\leq n$. Finally, we set $V_{\tau_i,s+1} = \{\gamma_i^\frown m: m>max(\ran(\gamma_i))\}$ for each $i<m$ and for every $i<m, e$, if $\tau_i$ has not been marked $e$-divergent, we set $U_e^{\tau_i} = U_e^{\tau_w}$ if $U_e^{\tau_w}$ has not been stopped, and $U_e^{\tau_i} = \emptyset$ otherwise. 

Let $G$ be $3$-generic. We want to find an enumeration $h:\omega \to G$ in $[f]$ which succeeds. We shall find such an enumeration in stages as $h = \bigcup \sigma_e$. Let $\sigma_0 = \epsilon$. Given $\sigma_e = f(\hat \tau)$, we pick $\sigma_{e+1}$ to ensure we beat $\varphi_e$ from computing a weak $3$-generic. Consider the $\Sigma^0_3$ set of strings $$X = \{str(f(\tau): f(\tau)\succ \sigma_e\text{ and } f(\tau) \text{ is marked $e$-divergent} \}.$$
If $G$ meets this set, so $\sigma \prec G$ has an enumeration in $f(\tau)\succ \sigma_e$ with $str(f(\tau))=\sigma$ which is marked $e$-divergent, we set $\sigma_{e+1} = f(\tau)$. In this case, all further extensions above $f(\tau)$ will be locking and so $\varphi_e^h$ is not total. 
On the other hand, if $G$ avoids this set, let $\sigma\prec G$ with $\gamma \succ \sigma$ implying $\gamma \not \in X$. Then we extend $\sigma_{e+1} = f(\tau)$ where $str(f(\tau)) = \sigma$. Corresponding to $f(\tau)$, the witnessing set $U_e^{\tau}$ never gets stopped (by hypothesis on $\sigma\prec G$ avoiding $X$). Therefore $U_e^{\tau}$ is a dense $\Sigma^0_3$ set of strings which $\varphi_e^h$ neither meets/avoids.
\end{proof}

This leaves open the natural question of whether this is tight.
\begin{query} Is there an $m<3$ and an $X\subseteq \omega$ such that $X$ is $m$-generic and every enumeration of $X$ computes a weak $3$-generic?
\end{query}

We end this section with an observation on the distribution of generics in the Turing degrees. Suppose $X$ is $n$-generic and $Y$ is $m$-generic with $n<m$. In \cite{belowgenerics}, it is noted that if $X \leq_T Y$ and $2\leq n$, then $X$ is $m$-generic as well. 

On the other hand, for any such $m$-generic set $Y$, there is a $1$-generic set $G>_T Y$, such that $G$ is of properly $1$-generic degree (no sets Turing equivalent to $G$ are $2$-generic). To see this, fix a $\Delta^0_2$ set $H$ which is $1$-generic. Then, $Y$ is $H$-generic, since $H\leq_T 0'$ and $m\geq 2$, and so $H\oplus Y$ will be $1$-generic. Moreover, $H\oplus Y$ is of properly $1$-generic degree since it does not form a minimal pair with $0'$. 

\section{Randomness}

Unlike with generics, we can show that enumerations fail to provide enough strength to compute randoms at a lower level. 
\begin{lem}\label{lem:w2r}
Every weakly $2$-random has an enumeration that does not compute a $1$-random.
\end{lem}
\begin{proof}
Let $X$ be weakly $2$-random. We will build an enumeration $f: \omega \to X$ which does not compute a $1$-random. Let $P_c$ be the $c$-incompressible strings, so $P_c = \{Z \subseteq \omega: \forall n K(Z\uh_n) > n-c\}$. If $Z$ is $1$-random, then $\exists c (Z \in P_c)$.

For every $e$ and $c$, we want to ensure $\varphi_e^f$ is partial or $\varphi_e^f \not \in P_c$. We will build $f$ in stages as $f = \bigcup_s \sigma_s$. 

At $s=0$, let $\sigma_0=\epsilon$. At stage $s+1$, if we have $s = \langle e,c \rangle$, we proceed as follows: Search for $n\in \omega, \tau \in X^{<\omega}$, with $\tau \succ \sigma_s$ such that $K(\varphi_e^\tau\uh_n)< n-c $. If such a $\tau$ exists, set $\sigma_{s+1}=\tau$ we are done. Otherwise consider:
$$ T = \{\rho: \exists \tau \in X^{<\omega} \; \tau \succ \sigma \land \rho \prec \varphi_e^\tau\}$$
Then $T \leq_e X$ and by the hypothesis, any path $Y \in [T]$ would be $1$-random. If $T$ had an isolated path, then that path would be total and $\leq_e X$, contradicting the fact that weak $2$-randoms are quasiminimal as shown by Cholak, Miller and Soskova \cite{randomquasi}. 

Therefore, $T$ has no isolated paths. Suppose $T$ had no deadends, then since $T \leq_e X \leq_e X\oplus \overline X$, we can find a perfect subtree $\hat T$ of $T$ such that $\hat T \leq_T X$. Barmpalias and Wang \cite{barmp} showed that no weakly $2$-random can compute such a perfect, pathwise random tree. Therefore $\hat T$ must have dead ends. Picking one of the dead ends forces $\varphi_e^R$ to be partial, and so we are done.

\end{proof}
Recall that the effective Hausdorff dimension of a set $X$ is defined as $$\operatorname{dim}_H(X) = \liminf_n \frac{K(X\uh_n)}{n}$$ where $K$ is the prefix free Kolmogorov complexity of a string. The authors in \cite{dimrand} show that dimension 1 sets are precisely the sets which differ from a $1$-random set on a density $0$ set of points. The following proof is implicit in their argument, and can even be shown directly with a little extra work.
\begin{lem}\label{lem:randomdim}
Every enumeration of a $1$-random computes a set of effective Hausdorff dimension $1$.
\end{lem}
\begin{proof}
If we have $n$ independent fair coin flips, with $X$ counting the number of heads, we have $\mu = n/2$ and $\sigma^2 = n/4$. Therefore by Chebeshev's inequality we get
$$ P(|X-\mu |\geq t) \leq \frac{n}{4t^2}$$
Setting $n= 2^k$ and $t=2^{2k/3}$, we therefore get 
$$  P(|X-\mu| \geq 2^{2k/3}) \leq \frac{2^{-k/3}}{4}$$
Let $I_s$ be the interval $[2^{s+1},2^{s+2}]$. If $X$ is $1$-random, there is a $k$, such that for every $s>k$, we have $2^s - 2^{2s/3} \leq |X \cap I_s| \leq 2^s  + 2^{2s/3}$.

We can define $Y = \bigcup_{s=k}^\infty \sigma_s$ as follows: Let $\sigma_k = X\uh_{2^{k+1}}$. For $s\geq k$, given $\sigma_s$, wait for the enumeration $f$ to have $2^{s} - 2^{\frac{2s}{3}}$ many $1$'s enter in the interval $I_s$, and extend $\sigma_{s+1} \succ \sigma_s$ to contain precisely these $1$'s. 

Clearly $Y$ is computable from $f$. $Y$ agrees with $X$ on a set of density $1$, and so is coarsely equivalent to a $1$-random. Therefore by \cite{dimrand}, we have $dim_{\mathcal H}(Y)=1$.
\end{proof}
\section{Large enumeration cones}
Recall that a large e-cone is the Turing degrees which enumerate a quasiminimal i.e.\ an enumeration cone not contained in a Turing cone. We can now answer Slaman and Sorbi's question affirmatively. A weak $2$-generic is quasiminimal but every total degree above it computes a weak $2$-generic and so cannot be of minimal Turing degree. We can show a stronger result on large e-cones.

\begin{lem} Let $A$ be $2$-generic relative to $B$. Then $A$ and $B$ form a $\Sigma_1$ minimal pair.
\end{lem}
\begin{proof}
Suppose $W = W_e^A = W_i^B$. We need to show that $W$ is c.e. Consider the $\Sigma^0_2(B)$ set
$$ X =\{\sigma : \exists x,s \forall t \; x \in W_{e,s}^\sigma \text{ and } x \not \in W_{i,t}^B \}.$$
By hypothesis, $A$ must avoid this set. Let $\sigma \prec A$ be such that $\forall \tau \succ \sigma \; \tau \not \in X$. Then $W$ is c.e.\ since $W = \{x : \exists \tau \succ \sigma \; x \in W_e^\tau\}$.
\end{proof}
\begin{cor} There is a large e-cone which does not contain an initial segment of the Turing degrees.
\end{cor}
\begin{proof}Let $X$ be $3$-generic and let $Y \geq_e X$ be total. Then $Y$ computes a $\Delta^0_3$ escaping function and hence a $2$-generic, the halves of which are mutually $2$-generic relative to each other. Then the halves form a minimal pair in the $e$-degrees and so cannot both lie in the cone above $X$.
\end{proof}
  
We can extend this analysis to any level of genericity by moving beyond escaping functions and exploiting the full power of enumerations.

\begin{defn} We say $X \subseteq \omega$ is thick if there is no $n$ such that $2n \not \in X$ and $2n+1 \not \in X$. A string $\sigma \in 2^{<\omega}$ is thick if it has even length and can be extended to a thick set.
\end{defn}
Let $g:2^{<\omega} \to \omega$ be a bijection such that if $|\sigma|<|\tau|$ then $g(\sigma)<g(\tau)$.
\begin{thm}\label{thm:fattogeneric}If $\mathcal F \subseteq  2^\omega$ is comeager and upwards closed under $\leq_T$, then $\mathcal F$ contains a large e-cone. 
\end{thm}
\begin{proof} Since $\mathcal F$ is comeager, $\exists X\subseteq \omega$ such that every $X$-generic set is in $\mathcal F$. We  shall build a thick set $Y$ of quasiminimal e-degree, such that every enumeration of $Y$ computes an $X$-generic. 

We build $Y$ in stages as $Y = \bigcup Y_s$. Let $Y_0 = \emptyset$. At stage $2s+1$ we check if there is a thick $\tau \succ Y_{2s}$ such that $2n, 2n+1 \in \Gamma_s^\tau$ for some $n$. If such a $\tau$ exists, we set $Y_{2s+1} = \tau$, ensuring that $\Gamma_s^Y$ is not total. Otherwise let $Y_{2s+1} = Y_{2s}$. 

At stage $2s+2$, we proceed as follows: let $n$ be the least such that $\forall \sigma \in 2^n \; g(\sigma)\times 2 >|Y_{2s+1}|$. Then there is a level $N>n$ such that each string $\sigma \in 2^n$ has an extension $\tau_\sigma \in 2^N$ which meets $W_s^X$ or avoids $W_s^X$. Then we define $Y_{2s+2}$ to be $01$ or $10$ at the bits $2g(\sigma^\frown \gamma), 2g(\sigma^\frown \gamma)+1$ depending on if $\sigma^\frown \gamma0 \prec \tau_\sigma$ or $\sigma^\frown \gamma1 \prec \tau_\sigma$. The values of $Y_{2s+2}$ at all other locations upto $max_{\tau \in 2^N} 2g(\tau)+1$ can be filled with $1$'s. This ends the construction.

Suppose $\Gamma_s^Y$ is total with $Y$ thick, say $A\oplus \overline A = \Gamma_s^Y$. Then $A$ is computable: $n \in A \iff \exists \tau \succ Y_{2s} \;$ with $\tau$ thick and $2n \in \Gamma_s^\tau$ and $n \not \in A \iff \exists \tau \succ Y_{2s}$ with $\tau$ thick and $2n+1 \in \Gamma_s^\tau$. (thickness is closed under unions)

Now given an enumeration $f$ of $Y$, we can compute an $X$ $1$-generic as follows. Staring at $X_0 = \emptyset$, we set $X_{s+1} = X_{s}^\frown i$ where $i$ is $0$ if we see a $1$ appear at position $g(X_s)*2+1$ before a $1$ appears in position $g(X_s)*2$ in the enumeration $f$ and $i=1$ otherwise.
\end{proof}

\begin{defn}
We say that $Y$ is thick $X$ $1$-generic if it meets or avoids every thick c.e.\ in $X$ set of strings.
\end{defn}
Note that the proof above shows that thick $1$-generics are quasiminimal. Moreover it is also thickly dense to have the property that there levels $n<N$ such that every string at level $n$ meets or avoids $W_s^X$ by level $N$ and that $\sigma$ codes these extensions as in the proof above. The set of such strings $\sigma$ are $c.e.$ in $X'$. Therefore taking $Y$ to be any thick $X'$-generic in the proof above would work.
\begin{lem}
Every enumeration of a weak $n$-thick generic computes a weak $n$-thick generic.
\end{lem}
\begin{proof}
If $W$ is a dense $\Sigma^0_n$ set of strings, then it is $\Delta^0_n(X)$ to find, for any level $m$, and $\forall \sigma \in 2^m$ a $\tau_\sigma \succ \sigma$ with $\tau_\sigma \in W$. It is thickly dense and $\Sigma^0_n$ to code the path to take these extensions as in Theorem~\ref{thm:fattogeneric}. Therefore decoding an enumeration as above will result in a fat generic.
\end{proof}

Under the embedding $\iota: X\to X \oplus \overline X$, a Turing degree can be quite high up in the e-degrees. It is interesting to note the following result of Jockush: Every Turing degree contains a semirecursive set $A$, and semirecursive sets (being halves of $\mathcal K$-pairs) are quasiminimal. Therefore every Turing degree contains a quasiminimal set, and being high up as a set in the Turing degrees does not prevent quasiminimality. 

In his thesis \cite{almosteveryquasi}, Lagemann showed that almost every set is quasminimal. This was strengthed by Cholak, Miller and Soskova \cite{randomquasi} to show that every weakly $2$-random is quasiminimal.

Kjos-Hanssen \cite{dncsubset} showed that computing an infinite subset of a $1$-random is equivalent to computing a DNC function. Therefore there is a large e-cone contained in the DNC degrees. Moreover, Lemma~\ref{lem:randomdim} showed that enumerating a random lets you compute a effective dimension 1 sequence.
\begin{cor} There is a large e-cone contained in the degrees of dimension $1$ sets.
\end{cor}
On the other hand, there are no large e-cones in the PA degrees.
\begin{defn} A $DNC_2$ multifunction is a map $F: \omega \to \{ \{0\},\{1\},\{0,1\}\}$ such that for every $n\in \omega$ we have $\varphi_n(n)\not \in F(n)$.
\end{defn}
\begin{lem}(J. Miller) If $F$ is a $DNC_2$ multifunction and $G_F \leq_e X$, then $X$ is not quasiminimal.
\end{lem}
\begin{proof} Suppose the graph of $F$ is enumeration below $X$. Let $g$ be the computable function defined by :
$$
\varphi_{g(e,i)}(n) = 
\begin{cases} 0 & e \in 0'\\
\varphi_i(n) & \text{otherwise}
\end{cases}
$$
Now there are two possibilities: 
\begin{enumerate}
\item $\exists e \not \in 0' \forall i F(g(e,i)) $ is single valued. In this case we would have a $DNC_2$ function, a total non-computable object, enumeration below $X$.
\item $\forall e \not \in 0' \exists i F(g(e,i))$ is multivalued. Then $e \not \in 0' \iff \exists i $ and two axioms witnessing $F$ is multivalued at $F(g(e,i))$ and so $\overline K \leq_e X$.
\end{enumerate}
\end{proof}
\begin{thm} Every large e-cone has a set which is not of PA degree.
\end{thm}
\begin{proof}
Let $X \subseteq \omega$ be quasiminimal. We build an enumeration $f:\omega \to X$ in stages as $f= \bigcup_s \sigma_s$.  For every $e$, we want to meet the requirements 
\[ \mathcal R_e : \exists n \; \varphi_e^f(n)\uparrow \text{ or } \varphi_e^f(n)\downarrow \not \in \{0,1\} \text{ or } \varphi_e^f(n) \downarrow = \varphi_n(n)\downarrow
\]. 
Let $\sigma_0 = \epsilon$. Given $\sigma_s$, we meet $\mathcal R_s$ at stage $s+1$ as follows: Check if $ \exists n, \exists \tau \in X^{<\omega}, \tau \succ \sigma_s$ such that $\forall \gamma \in X^{<\omega}, \gamma \succ \tau$ we have $\varphi_e^\gamma(n)\uparrow$ or $\varphi_e^\gamma(n)\downarrow \not \in \{0,1\}$. If there is such a $\tau$, let $\sigma_{s+1} = \tau$. Otherwise consider the tree $$T = \{ \varphi_e^\gamma : \gamma \succ \sigma_s, \gamma \in X^{<\omega}\}.$$ By the hypothesis, any path through $T$ is $DNC_2$ and $T$ has no deadends. Moreover $T\leq_e X$ and so if $T$ has an isolated path, that path would be $\leq_e X$, a contradiction to the quasiminimality of $X$. Hence $T$ is a perfect tree with no deadends. Using $T$ we get a $DNC_2$ multifunction $F$ such that $G_F \leq_e X$ contradicting the lemma above.

Therefore the tree must either have a dead end allowing us to force partiality or must agree with the diagonal at some point allowing us to force non-$DNC_2$, and we would have successfully choosen $\sigma_{s+1}$ accordingly. 
\end{proof}

\begin{lem}
There is a large e-cone not containing any set of $1$-generic Turing degree and a large e-cone not containing any set of $2$-random degree.
\end{lem}
\begin{proof}
Let $X$ be weak $2$-random. Then $X$ is quasiminimal by \cite{randomquasi} and all the total degrees above it compute a DNC function, and hence cannot be $1$-generic. 

Let $Y$ be weak $2$-generic, then $Y$ is quasiminimal and every enumeration of $Y$ computes a weak $2$-generic. Therefore, no such enumeration can be of $2$-random degree since $2$-randoms form minimal pairs with weak $2$-generics as noted in \cite{minpair}. 
\end{proof}


\section*{Acknowledgements}
I would like to thank my advisor, Joe Miller, for many helpful discussions and for suggesting the questions in Section~2. I also thank Mariya Soskova for bringing the question of Slaman and Sorbi to my attention. The only use of AI tools was in the proof of Lemma~\ref{lem:w2r} and with a part of the combinatorial argument in the proof of Theorem~\ref{thm:enumnosuf}.

\bibliographystyle{plain}
\bibliography{ref}
\end{document}